\RequirePackage{snapshot}
\documentclass[times]{elsarticle}

\biboptions{sort&compress}
\usepackage{jcomp}
\usepackage{framed,multirow}

\usepackage{amssymb}
\usepackage{latexsym}
\usepackage{url}
\usepackage{xcolor}
\definecolor{newcolor}{rgb}{.8,.349,.1}

\usepackage[T1]{fontenc}
\usepackage[utf8]{inputenc}
\usepackage{siunitx}
\usepackage{amsmath}
\usepackage{scalerel,stackengine}
\usepackage{graphicx}
\usepackage{tikz}
\usepackage{afterpage}
\usepackage{subcaption}
\usepackage{tabularx}
\usepackage{multirow}

\usepackage[normalem]{ulem}

\usepackage[
  vlined,
]{algorithm2e}

\renewcommand{\vec}[1]{\boldsymbol{#1}}
\newcommand{\mat}[1]{\mathbf{#1}}
\newcommand\equalhat{\mathrel{\stackon[1.5pt]{=}{\stretchto{%
    \scalerel*[\widthof{=}]{\wedge}{\rule{1ex}{3ex}}}{0.5ex}}}}
\journal{Journal of Computational Physics}

\begin{document}

\begin{frontmatter}

\title{A Variational Integrator for the Discrete Element Method}

\author[1]{David N. {De Klerk}\corref{cor1}}
\cortext[cor1]{Corresponding author:
  Email: \href{mailto:david.deklerk@glasgow.ac.uk}{david.deklerk@glasgow.ac.uk}}
\author[1]{Thomas {Shire}}
\author[1]{Zhiwei {Gao}}
\author[1]{Andrew T. {McBride}}
\author[1]{Christopher J. {Pearce}}
\author[1,2]{Paul {Steinmann}}

\address[1]{Glasgow Computational Engineering Centre, James Watt School of Engineering, University of Glasgow, Glasgow G12 8QQ, United Kingdom}
\address[2]{Institute of Applied Mechanics (LTM), Friedrich-Alexander Universität Erlangen-Nürnberg (FAU), Erlangen, Germany}

\begin{abstract}
  A novel implicit integration scheme for the Discrete Element Method (DEM)
  based on the variational integrator approach is presented.
  The numerical solver provides a fully dynamical description that, notably, reduces to an energy
  minimisation scheme in the quasi-static limit.
  A detailed derivation of the numerical method is presented for the Hookean contact
  model and tested against an established open source DEM package that uses the
  velocity-Verlet integration scheme.
  These tests compare results for a single collision, long-term stability
  and statistical quantities of ensembles of particles.
  Numerically, the proposed integration method demonstrates equivalent accuracy to
  the velocity-Verlet method.
\end{abstract}


\end{frontmatter}

\section{Introduction}

\begin{figure}
  \centering
  \includegraphics{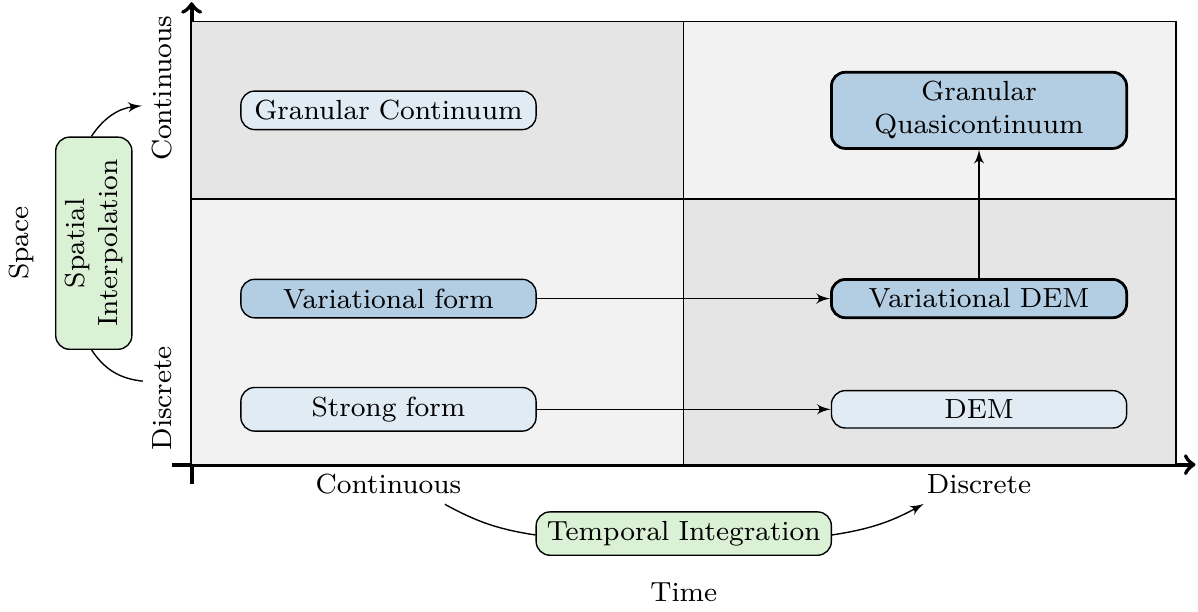}
  \caption{
    A classification of descriptions of granular material based on the treatment (continuous or discrete)
    of the temporal and spatial dimensions.
    The focus of this work is a variational integrator for the Discrete Element Method
    which will provide the appropriate mathematical setting for a granular Quasicontinuum Method.
  }
  \label{fig:outline}
\end{figure}

Various descriptions of granular materials are compared  in Figure~\ref{fig:outline}, where
they are classified by the treatment of the temporal and spatial dimensions, which can be
either continuous or discrete.
In the underlying Newtonian picture (strong form), the discrete spatial degrees of freedom (particle
position and orientation) are described by continuous functions of time which are
the solutions to Newton's second law.
The Discrete Element Method (DEM) \cite{Cundall1979}
-- a widely-adopted particle-level approach for simulating granular materials --
calculates the resultant force acting on each particle during distinct time steps resulting
from the discretisation of the time domain and solves
the governing equation of motion.
A continuum description (granular continuum) of spatially discrete systems is achieved via a
micro-to-macro transition.
For example Babic~\cite{Babic1997} proposed a coarse-graining method and derived a balance equation
that relates continuous functions of position to each other.
In practice, the micro-to-macro transitions for granular systems are often performed on
discrete-time/discrete-space DEM data~\cite{Miehe2004}, but in principle this can be achieved for the Newtonian
description too.

Unlike computational models of fluid dynamics or continuum mechanics, numerical
simulation of granular materials have not been able take advantage of developments in
spatial continuum modelling.
Granular materials display a variety of behaviours
which is often compared to the solid, fluid and gaseous phases of matter.
The solid-like phase is characterised by static packing and jamming,
the energetic gaseous state by pairwise collisions between particles,
and the intermediate fluid-like state by dense flows
\cite{Jaeger1996}.
Given the complexity and diversity of physical phenomena present in
granular materials, finding a universal continuum description for
granular material remains an open research question.
A local continuum description for dense granular
flow has been proposed \cite{GDRMidi2004, daCruz2005, Jop2006}
and have shown to be applicable in a range of situations.
However, this rheology has limitations and
fails to reproduce important non-local phenomena such as shear banding and arching
\cite{Jop2015, Kamrin2019}.
The former occurs a granular assembly is subjected to shear loading.
While there is significant particle rotation and relative motion inside the shear band
the remaining part of the assembly typically moves like a
rigid body \cite{Gao2013}. %

\begin{figure}
  \centering
  \begin{subfigure}[b]{0.25\textwidth}
    \centering
    \includegraphics[width=.8\textwidth]{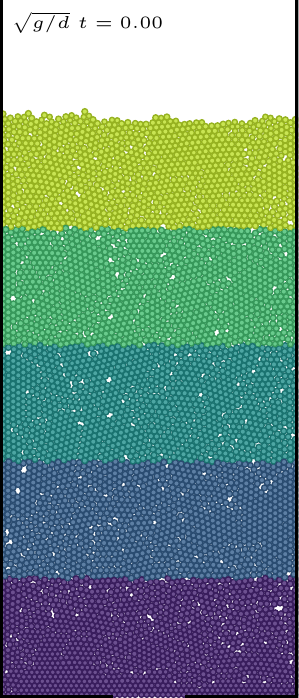}
  \end{subfigure}%
  \begin{subfigure}[b]{0.25\textwidth}
    \centering
    \includegraphics[width=.8\textwidth]{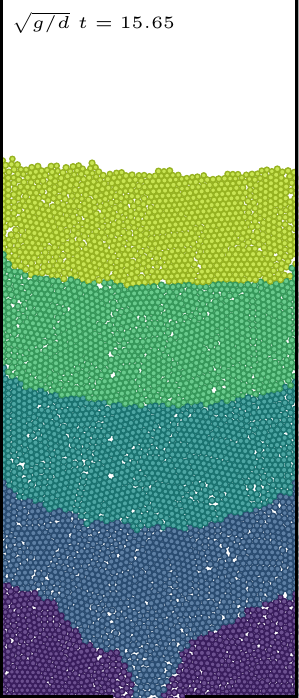}
  \end{subfigure}%
  \begin{subfigure}[b]{0.25\textwidth}
    \centering
    \includegraphics[width=.8\textwidth]{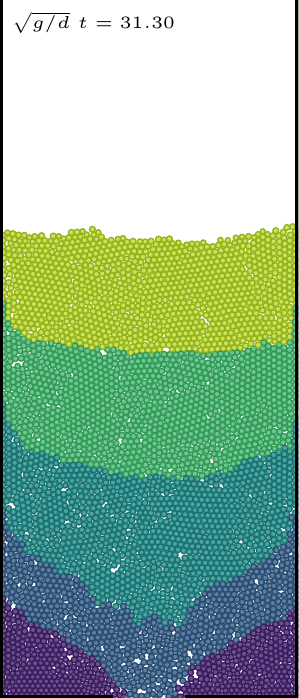}
  \end{subfigure}%
  \begin{subfigure}[b]{0.25\textwidth}
    \centering
    \includegraphics[width=.8\textwidth]{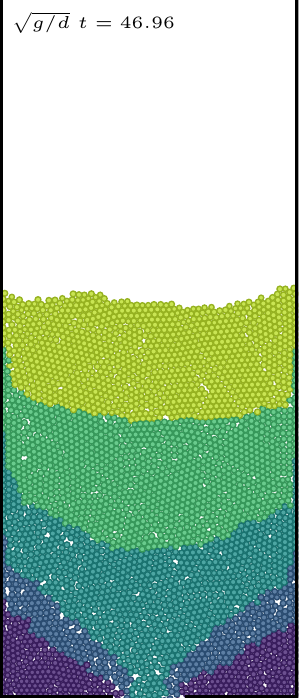}
  \end{subfigure}
  \caption{
    A DEM simulation of a draining silo.
    Many particles experience only small relative displacements for the
    largest part of the simulation, which is common in many DEM
    applications.
  }
  \label{fig:weinhart2016}
\end{figure}

In the absence of a complete continuum theory, many studies of granular
materials rely on discrete, particle level numerical simulations.
However, in the static or slow moving phase of granular materials,
a large number of particles may remain nearly stationary or behave in
a manner that could be described by a continuum model.
For instance, Figure~\ref{fig:weinhart2016}  shows a draining silo with
particles coloured by their initial vertical position.
Even at an advanced state of drainage, particles in certain regions still
approximately maintain their positions relative to their initial neighbours.

The Quasicontinuum (QC) method is a multiscale discrete-continuum method which allows for
a fully-resolved particle simulation where required,  and a more
efficient continuum description of material behaviour elsewhere.
Simulations are carried out in a continuous spatial domain and the method thus fits to the
top right quadrant of the diagram in Figure~\ref{fig:outline}.
The method was initially developed for crystalline atomistic simulations
\cite{Tadmor1996, Knap2001, Miller2002, Tadmor2005},
where  the arrangement of atoms is calculated so as to minimise the global potential energy
of the system using a suitable numerical technique, such as
  iterative energy minimisation methods.
Atoms exist throughout the domain, but the computational cost is reduced by two key features.
First, a series of representative atoms, or rep-atoms, are identified.
The density of rep-atoms is highest in regions of specific interest and gradually diffuses toward
regions of less interest.
Second, the energy density is estimated by so-called summation rules in regions bordered by
rep-atoms.
The displacement of non rep-atoms are updated by interpolating their positions
between rep-atoms.
In situations where the majority of atoms fall in regions of low interest, the degrees
of freedom of the simulation is greatly reduced which leads to improved simulation run time.

The objective here is to provide a temporal discretisation framework for the application of
the QC method to granular systems. However,
several challenges exist before a granular QC method can be realised.
For the most part, with exceptions, such as \cite{Kochmann2014, Amelang2015},
only quasi-static configurations are simulated in the QC framework and dynamics are not accounted for.
The original QC method was developed for quasi-static crystalline atomistic simulations that
minimises the inter-atomic potential energy of the system.
In the context of granular materials, a quasi-static simulation would restrict the
method's application to the solid-like state.
To recover the dynamics, and to stay consistent with the QC approach,
a novel integration method for the Discrete Element Method has been developed.
The method follows
Hamilton's principle in seeking the stationary point of the action.
The other hallmark of the QC method -- an efficient summation rule
-- will be addressed in future work.

In a time continuous setting, Hamilton's principle provides a variational scheme where the differential
equations governing a dynamical system can be derived by finding the trajectory that is the
stationary point of the action.
The Lagrange-d'Alembert principle is a generalisation of Hamilton's principle to non-holonomic
systems, and is therefore applicable here due to the dissipative nature of granular materials.
The classification in Figure~\ref{fig:outline},
identifies the Hamiltonian approach (variational) to be in the same category as the Newtonian one (strong form).
Variational integrators \cite{Marsden1999, Kane2000, Marsden2001, Lew2004a, Lew2004}
are a class of algorithms where the time continuous variational principles are discretised
to obtain time-stepping schemes for dynamical systems.
As a result, many of the important properties of Lagrangian mechanics carry over to these algorithms.
For instance, variational integrators conserve the generalised momentum of a system as a consequence of a
discrete version of Noether's theorem.
Of particular interest is an implicit integration scheme, outlined in \cite{Marsden2001},
that follows directly from Hamilton's principle in a  discrete setting.
When the quasi-static approximation is taken, i.e. by neglecting inertia, the method
simplifies to minimising the potential energy of the system -- precisely what is done  in
the atomistic simulations that inspired the QC method.
A variational integrator for DEM is the time discrete analogue to the variational format,
i.e. the time discrete description in the bottom right quadrant of Figure~\ref{fig:outline},
and provides a way to proceed towards the top right quadrant.

While proposing variational integrators per se is not new, the bespoke  application to the Discrete
Element Method is novel. This is a crucial step towards a granular Quasicontinuum method.
To achieve this objective, a benchmark against current DEM solvers is needed before
addressing the other challenges mentioned above.
The remainder of the paper is structured as follows.
Section 2 provides a detailed derivation of the variational integrator for dissipative
systems and extends its application to the Hookean contact model in DEM.
Section 3 discusses the implementation of the solver, shows results of numerical
experiments and comparisons with established DEM codes.
Section 4 is dedicated to the final discussion and conclusions.

\section{Numerical Integration}

The velocity-Verlet method \cite{Verlet1967}
is popular in molecular dynamics and is also widely used in DEM.
The same method is also known as the Str{\"o}mer method and the leapfrog method, depending on the context where it is used \cite{Hairer2003}.
It has been shown that the velocity-Verlet method and many of its variants can be derived using the variational
integrator approach \cite{Ruth1983, Leimkuhler1994, Hairer2003} and therefore
inherits the properties of variational integrators mentioned in the introduction.
However, the velocity-Verlet method is explicit and tailored toward solving Netwon's
equations in the strong form (see Figure~\ref{fig:outline}) and, as discussed above,
a variational approach is preferred for the Quasicontinuum method.

In practice, DEM simulations are carried out over time periods many orders of magnitude
larger than the duration of a single contact which leads to a trade off between the duration
and the accuracy or stability of the simulation.
To ensure accurate particle trajectories, an integration time step needs to be selected that
is much smaller than the duration of a contact.
Choosing the optimal integration time step has been the topic of substantial
research \cite{OSullivan2004, Washino2016, Otsubo2017a}.
Implicit integration schemes have been proposed for DEM (see for instance \cite{Ke1995, Samiei2013}).
However the same limitation on the maximum time step applies, and with the added computational cost
of implicit schemes these methods typically results in longer simulation times
than explicit schemes.
An approach to solve DEM by minimising the potential energy was proposed in
\cite{Krijgsman2016}, however this method was restricted to quasi-static configurations.

\subsection{Variational Integrators}

\begin{figure}
  \centering
  \begin{subfigure}[b]{0.4\textwidth}
    \centering
    \includegraphics{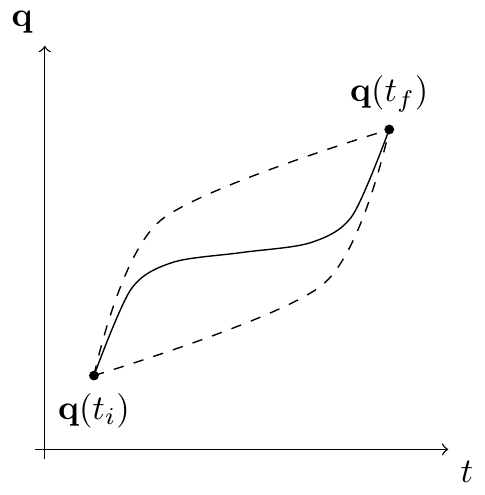}
    \caption{}
  \end{subfigure}
  \begin{subfigure}[b]{0.4\textwidth}
    \centering
    \includegraphics{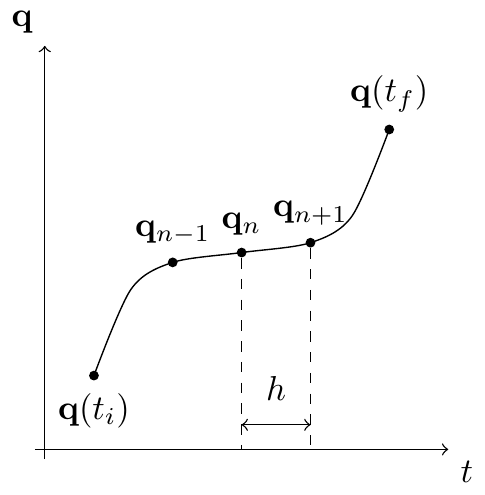}
    \caption{}
  \end{subfigure}
  \caption{Hamilton's Principle of least action (a) is a continuous method for solving
    the trajectory of a dynamical system by finding the path that is the stationary point
    of the action.
    Variational integrators (b) are the discrete realisation of this principle.}
  \label{fig:hamilton}
\end{figure}

The numerical integration scheme presented here follows Kane et al. \cite{Kane2000}.
The Lagrange-d'Alembert principle (see Figure~\ref{fig:hamilton} (a)) is used to derive a
second-order accurate integrator
for the equations of motion of a general dynamical system.
The continuous formulation of this principle states that for a system
under the influence of a non-conservative generalised force
$\vec{Q}(\vec{q}, \dot{\vec{q}})$, the sum of the variation of the action
($S = \int L\ dt$)
and the total work performed by the non-conservative forces is zero, that is
\begin{equation}
  \underbrace{\delta \int_{t_i}^{t_f} L(\vec{q}, \dot{\vec{q}}) \ dt}_{ \delta S}
  + \int_{t_i}^{t_f} \vec{Q}(\vec{q}, \dot{\vec{q}}) \cdot \delta \vec{q} \ dt = 0 .
  \label{eq:dalambert}
\end{equation}
Here $\vec{q}$ and $\dot{\vec{q}} = d\vec{q}/dt$ are the generalised coordinates
and velocities, respectively,
and the Lagrangian is given by $L(\vec{q}, \dot{\vec{q}}) = T(\dot{\vec{q}}) - V(\vec{q})$,
where $T(\dot{\vec{q}})$ and $V(\vec{q})$ are the kinetic and potential energy of the system,
respectively.

A generalised coordinate can be any parameter that specifies the configuration
of the system. For discrete particles these are the coordinates and angles
that specify their position and orientation.
The corresponding generalised forces are forces and torques.

In the absence of non-conservative forces ($\vec{Q}=\vec{0}$), the second term in Eq.
(\ref{eq:dalambert}) is zero and the Lagrange-d'Alembert principle is equivalent to
Hamilton's principle of least action.
The Lagrange-d'Alembert principle will be required to formulate an integrator for
DEM, because of the dissipative terms in the contact model.
Since Hamilton's principle is a special case, we will refer to it
in the following discussion, when appropriate.

In order to find the trajectory that a system will follow in the time continuous case,
the calculus of variations is used to find the stationary point of the action.
For a time discrete formulation,
the trajectory  $\vec{q}(t)$ is decomposed into $N$ time steps of
length $h$ and labelled $\{\vec{q}_0, \cdots, \vec{q}_N\}$ as depicted in Figure~\ref{fig:hamilton} (b).
A discrete Lagrangian is defined as the numerical approximation
of the integral over the time step and is given by,
\begin{equation}
  L^\alpha_d(\vec{q}_k, \vec{q}_{k+1}, h)
  = \int_t^{t+h} L(\vec{q}, \dot{\vec{q}})\ dt
  \approx h L\left(\vec{q}_{k+\alpha}, \frac{\vec{q}_{k+1} - \vec{q}_k}{h}\right) ,
  \label{eq:dislag}
\end{equation}
where
$\vec{q}_{k+\alpha} = [1-\alpha]\vec{q}_k + \alpha\vec{q}_{k+1}$.
The parameter $\alpha$ is often chosen as $0$ or $1/2$ which correspond to
the left hand rule or midpoint rule, respectively.
The former leads to a first-order accurate integrator and the latter increases the
accuracy to second-order.
The discrete action is the sum over the $N$ time steps,
\begin{equation}
  S_d  = \sum_{k=0}^{N-1}  L^\alpha_d(\vec{q}_k, \vec{q}_{k+1}, h).
\end{equation}

The discrete Lagrange-d'Alembert principle \cite{Kane2000} is given by
\begin{equation}
  \delta \sum_{k=0}^{N-1} L^\alpha_d(\vec{q}_{k}, \vec{q}_{k+1}) +
  \sum_{k=0}^{N-1}
  \Big[
    \vec{Q}^-(\vec{q}_{k}, \vec{q}_{k+1})\cdot \delta \vec{q}_{k} +
    \vec{Q}^+(\vec{q}_{k}, \vec{q}_{k+1})\cdot \delta \vec{q}_{k+1}
  \Big]
  = 0
  \label{eq:disdalambert}
\end{equation}
where,
\begin{align}
  \vec{Q}^-_d(\vec{q}_{k}, \vec{q}_{k+1})
  &= \frac{h}{2} \vec{Q}\left(\vec{q}_{k+\alpha}, \frac{\vec{q}_{k+1} - \vec{q}_{k}}{h}\right),
  \\
  \vec{Q}_d^+(\vec{q}_{k}, \vec{q}_{k+1})
  &= \frac{h}{2} \vec{Q}\left(\vec{q}_{k+1-\alpha}, \frac{\vec{q}_{k+1}-\vec{q}_{k}}{h}\right).
\end{align}

The equivalent of the Euler-Lagrange equations can be derived using
Hamilton's principle of stationary action.
The dynamics of the system will ensure that the variation in the action,
$\delta S_d$, remains  zero for independent variations in $\delta \vec{q}_k$
and $\delta \vec{q}_{k+1}$, that is
\begin{equation}
  \delta S_d =
  \sum_{k=0}^{N-1}
  \delta \vec{q}_k\cdot \frac{d}{d\vec{q}_k}L^\alpha_d(\vec{q}_k, \vec{q}_{k+1}, h)\ +
  \sum_{k=0}^{N-1}
  \delta \vec{q}_{k+1}\cdot \frac{d}{d\vec{q}_{k+1}}L^\alpha_d(\vec{q}_k, \vec{q}_{k+1}, h) .
  \label{eq:hamilton1}
\end{equation}
To prevent confusion with derivatives, new notation is introduced
such that $D_1L^\alpha_d$ and $D_2L^\alpha_d$ are the derivative of the first and second argument of $L^\alpha_d$,
respectively.
Then, the index for the sum in the second term is changed to $k+1 \rightarrow k$.
The expression for Eq. (\ref{eq:hamilton1}) now becomes,
\begin{equation}
  \delta S_d =
  \sum_{k=0}^{N-1}
  \delta \vec{q}_k \cdot D_1 L^\alpha_d(\vec{q}_k, \vec{q}_{k+1}, h)\ +
  \sum_{k=1}^{N}
  \delta \vec{q}_{k}\cdot  D_2 L^\alpha_d(\vec{q}_{k-1}, \vec{q}_k, h) .
\end{equation}
Since $\delta \vec{q}_0 = \delta \vec{q}_N = 0$, the first sum can start
at $k=1$ and the second can be terminated at $N-1$.
Now, since both summations are carried out over the same range, the expression can be
factorised, as
\begin{equation}
  \delta S_d =
  \sum_{k=1}^{N-1}
  \delta \vec{q}_k
  \cdot
  \Big[
  D_1 L^\alpha_d(\vec{q}_k, \vec{q}_{k+1}, h)  +
  D_2 L^\alpha_d(\vec{q}_{k-1}, \vec{q}_k, h)
  \Big] .
  \label{eq:ds}
\end{equation}

The condition $\delta S_d = 0$ can be enforced by requiring that the term in the brackets
be zero, that is,
\begin{equation}
  D_1 L^\alpha_d(\vec{q}_k, \vec{q}_{k+1}, h) + D_2 L^\alpha_d(\vec{q}_{k-1}, \vec{q}_k, h) = \vec{0}  , \label{eq:del}
\end{equation}
which is the discrete form of the Euler-Lagrange equation.

In the general case when dissipative forces are present (i.e. $\vec{Q} \neq \vec{0}$), the
second term in Eq. (\ref{eq:disdalambert}) can be manipulated using the same steps as above to
obtain,
\begin{equation}
  \sum_{k=0}^{N-1}
  \Big[
    \vec{Q}_d^-(\vec{q}_{k}, \vec{q}_{k+1}, h) \cdot \delta \vec{q}_{k} +
    \vec{Q}_d^+(\vec{q}_{k}, \vec{q}_{k+1}) \cdot \delta \vec{q}_{k+1}
  \Big]
  =
  \sum_{k=1}^{N-1}  \delta \vec{q}_{k} \cdot
  \Big[
    \vec{Q}_d^-(\vec{q}_k, \vec{q}_{k+1}, h)  +
    \vec{Q}_d^+(\vec{q}_{k-1}, \vec{q}_{k})
  \Big] \ .
\end{equation}
The sum over $k$ and $\delta \vec{q}_k$ can be factored with the terms in (\ref{eq:ds}), which leads to the discrete Euler-Lagrange equation,
\begin{equation}
  D_1 L^\alpha_d(\vec{q}_k, \vec{q}_{k+1}, h) + D_2 L^\alpha_d(\vec{q}_{k-1}, \vec{q}_k, h)
  +
  \vec{Q}_d^-(\vec{q}_k, \vec{q}_{k+1}, h)  + \vec{Q}_d^+(\vec{q}_{k-1}, \vec{q}_{k}, h) = \vec{0}  .
  \label{eq:dfel}
\end{equation}

Both Eqs. (\ref{eq:del}) and (\ref{eq:dfel}) are second-order 
equations, but
a system of two first-order equations can be constructed by introducing the generalised
momentum.
The momentum in the time continuous case is defined by,
$\vec{p}(t) = \partial L/\partial \dot{\vec{q}}$.
Similarly in the time discrete setting, the momentum at step $k$ is given by
\begin{equation}
  \vec{p}_k \equiv D_2 L^\alpha_d(\vec{q}_{k-1}, \vec{q}_k, h) .
  \label{eq:momdef}
\end{equation}
The momentum can be used to evaluate $\vec{Q}^+$ at time step $k$,
\begin{equation}
  \vec{Q}_d^p(\vec{q}_{k}, \vec{p}_{k}) = \frac{h}{2} \vec{Q}\left(\vec{q}_{k}, \frac{\vec{p}_k}{m}\right) = \vec{Q}^+(\vec{q}_{k-1}, \vec{q}_{k}, h) .
\end{equation}
The first update equation is obtained by
substituting  the
definition for the momentum (\ref{eq:momdef}) at step $k$,
into the discrete Euler-Lagrange equation (\ref{eq:dfel}),
and the second is the expression for the momentum at step $k+1$.
The pair of first-order update equations is given by,
\begin{align}
  \vec{R}(\vec{q}_k, \vec{q}_{k+1}, \vec{p}_k, h)
  &\equiv
  \vec{p}_k
  +
    D_1 L^\alpha_d(\vec{q}_{k}, \vec{q}_{k+1}, h) +
    \vec{Q}_d^-(\vec{q}_k, \vec{q}_{k+1}, h)  +
    \vec{Q}_d^p(\vec{q}_{k}, \vec{p}_{k})
  =
  \vec{0},
  \label{eq:update:q} \\
  \vec{p}_{k+1} &= D_2 L^\alpha_d(\vec{q}_k, \vec{q}_{k+1}, h). \label{eq:update:p}
\end{align}
The update scheme, $(\vec{q}_k, \vec{p}_k) \mapsto (\vec{q}_{k+1}, \vec{p}_{k+1})$, requires that
the new position, $\vec{q}_{k+1}$ be calculated using an implicit scheme in Eq.
(\ref{eq:update:q})
and then explicitly calculating the new momentum $\vec{p}_{k+1}$ using (\ref{eq:update:p}).

The implicit scheme for updating $\vec{q}$ is obtained by expanding
(\ref{eq:update:q}) around $\vec{q}_{k+1}$,
\begin{equation}
  \vec{R}(\vec{q}_k, \vec{q}_{k+1}^n, \vec{p}_k, h) + \mat{K}(\vec{q}_k, \vec{q}_{k+1}^n, \vec{p}_k, h) \Delta \vec{q}_{k+1}^n = \vec{0},
  \label{eq:implicit:step}
\end{equation}
where $\vec{q}^n_{k+1}$ is the previous estimate for $\vec{q}_{k+1}$
and
$\Delta \vec{q}_{k+1}^n = \vec{q}_{k+1}^{n+1} - \vec{q}_{k+1}^n$
is the change required to improve the estimate.
The stiffness is given by,
\begin{align}
  \mat{K}(\vec{q}_k, \vec{q}_{k+1}, h) &=
  \frac{\partial}{d\vec{q}_{k+1}} \vec{R}(\vec{q}_k, \vec{q}_{k+1}, h) \\
  &=
  \frac{\partial^2}{\partial \vec{q}_k \partial \vec{q}_{k+1}} L^\alpha_d(\vec{q}_k, \vec{q}_{k+1}, h) +
  \frac{\partial}{\partial \vec{q}_{k+1}} \vec{Q}^-_d(\vec{q}_k, \vec{q}_{k+1}, h)
\end{align}
The initial guess for $\vec{q}_{k+1}$ can be estimated using the momentum at the $k$-th time
step, $\vec{q}^0_{k+1} = \vec{q}_k + h \vec{p}_k/m$.

\subsection{Integration for DEM}

DEM is characterised by treating particles as rigid bodies with `soft' contacts
where overlap between particles is allowed and
inter-particle forces are expressed as a function of the overlap (denoted by $\delta_{ij}$).
Different contact models have been proposed (for instance see \cite{Rojek2018} for a recent review),
but for simplicity and without loss of generality the Hookean contact model \cite{Silbert2001c, Rojek2018}
is adopted.
Specifically, the normal and tangential forces between particles, expressed in the global
coordinate system, are calculated using
\begin{align}
  \vec{F}_{n_{ij}} &=
  k_n \delta_{ij} \vec{n}_{ij} - \gamma_n m_\mathrm{eff} \vec{v}_{n_{ij}} ,
  \label{eq:dem_contact_normal}
  \\
  \vec{F}_{t_{ij}} &=
  k_t \delta_{ij} \vec{t}_{ij} - \gamma_t m_\mathrm{eff} \vec{v}_{t_{ij}},
  \label{eq:dem_contact_tan}
\end{align}
where
$k_n$ and $k_t$ are the normal and tangential spring stiffness,
$\gamma_n$ and $\gamma_t$ are the normal and tangential damping coefficients and
$m_\mathrm{eff} = m_i m_j / [m_i + m_j]$ is the effective mass of the contact.
The overlap, normal and tangential components of the velocity, are given by,
\begin{align}
  \vec{v}_{n_{ij}} &= \left[\vec{v}_{ij}\cdot \vec{n}_{ij}\right] \vec{n}_{ij} ,\\
  \vec{v}_{t_{ij}} &= \vec{v}_{ij} - \vec{v}_{n_{ij}} - \frac{1}{2}
  \left[\vec{\omega}_i + \vec{\omega}_j\right]\times \vec{r}_{ij} ,
\end{align}
respectively,
where
$\vec{r}_{ij} = \vec{r}_i - \vec{r}_j$ is the relative position of the particles,
$\vec{n}_{ij} = \vec{r}_{ij}/|\vec{r}_{ij}|$ is the unit vector normal to the contact,
$\vec{t}_{ij}$ is unit vector tangential to the contact,
$\delta_{ij} = d - |\vec{r}_{ij}|$ is the overlap between the particles with diameter $d$,
$\vec{v}_{ij} = \vec{v}_i - \vec{v}_j$ is the relative velocity
and $\vec{\omega}$ is the angular velocity.

The dynamics of particle $i$ is governed by the resultant force and torque,
\begin{align}
  \vec{F}_i &= \vec{F}_i^\mathrm{ext} + \sum_j\left[\vec{F}_{n_{ij}} + \vec{F}_{t_{ij}}\right], \\
  \vec{\tau}_i &= -\frac{1}{2} \sum_j\left[\vec{r}_{ij}\times \vec{F}_{t_{ij}}\right],
\end{align}
where
$\vec{F}_i^\mathrm{ext}$ are any external forces on particle $i$
and
the sum $j$ is carried out over all particles that are in contact with
$i$, i.e. for which $\delta_{ij} > 0$.

The DEM method can be cast into the Lagrangian formulation where the Lagrangian for a system of
$N_p$ discrete particles is given by,
\begin{equation}
  L(\vec{q}, \dot{\vec{q}}) =
  \frac{1}{2}\dot{\vec{q}}^T \mat{M} \dot{\vec{q}}
  -
  V(\vec{q}) ,
  \label{eq:dem:lagrangian}
\end{equation}
where $\vec{q}$ is a $6N_p$ real column vector that represents the degrees of freedom
of all $N_p$ particles.
The integrator needs to account for the position and orientation
of each particle,
so a reasonable choice is to group the vector in rows of $6$, where the
first $3$ entries and last $3$ entries represent the position and
angular degrees of freedom, respectively.
As before, the generalised velocity is $\dot{\vec{q}} = d\vec{q}/dt$ which is
therefore composed of the linear and angular velocity.
The generalised momentum,
$\vec{p}$, contains both the linear and angular
momentum.
A component of $\vec{p}$ is given by $p =\partial L/\partial \dot{q}$,
where $\dot{q}$ is the corresponding component of $\dot{\vec{q}}$.
The mass matrix $\mat{M}$ is a
$6N_p\times6N_p$ diagonal matrix with blocks $\mat{M}_i = \mathrm{diag}([m_i\ m_i\ m_i\ I_i\ I_i\ I_i])$.
Here  $m_i$ and $I_i$ are the particle mass and moment of inertia, respectively, of
particle $i$.

\begin{figure}
  \centering
  \includegraphics{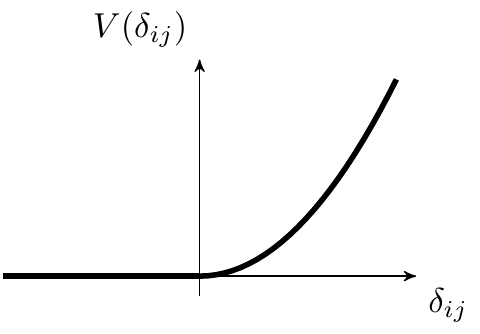}
  \caption{The Hookean inter-particle potential function as a function of particle overlap.}
  \label{fig:vij}
\end{figure}

The first term in Eq. (\ref{eq:dem:lagrangian}) accounts for the total kinetic energy of the
system.
The potential energy due to a Hookean contact between particles $i$ and $j$ is given by
\begin{align}
  V_{ij} =
  \begin{cases}
  \frac{k_n}{2} \left[ \delta_{ij}\right]^2, & \text{if } \delta_{ij} > 0 \\
  0, &\text{otherwise}.
  \end{cases}
\end{align}
The potential function is illustrated in Figure~\ref{fig:vij}.
This formulation can be expanded to other contact models, for instance
a Hertz-Mindlin contact model can be implemented by using
$V_{ij} \propto 2/5\ [\delta_{ij}]^{5/2}$ for $\delta_{ij}>0$.
The potential energy for the entire system is the sum of the potentials over all
the particles in contact with each other
and, assuming a gravitational acceleration $\vec{g} \equalhat [0, 0, -g]$, the gravitational potential energy
$V = \sum_i \sum_{j<i} V_{ij} + \sum_i m_igz_i$.
The generalised non-conservative forces $\vec{Q}$ in DEM are friction and
velocity dependent damping terms in Eqs.~(\ref{eq:dem_contact_normal}-\ref{eq:dem_contact_tan}).

Following the prescription in Eq. (\ref{eq:dislag}), the discrete Lagrangian for DEM is given by,
\begin{equation}
  L^\alpha_d(\vec{q}_{k}, \vec{q}_{k+1}, h) =
  \frac{1}{2h}
  \left[\vec{q}_{k+1} - \vec{q}_{k}\right]^T
  \mat{M}
  \left[\vec{q}_{k+1} - \vec{q}_{k}\right]
  -
  h
  V(\vec{q}_{k+\alpha})
  .
\end{equation}
This allows one to simplify the update scheme for the integrator. The vector term in Eq.~(\ref{eq:implicit:step}) becomes,
\begin{equation}
    \vec{R}(\vec{q}_k, \vec{q}_{k+1}^n, \vec{p_k}, h)
    =
    \vec{p}_k
    -
    \frac{1}{h}\mat{M}\left[\vec{q}^{n}_{k+1} - \vec{q}_{k}\right]
    -
    h
    \left[
    [1-\alpha]
    \frac{\partial V}{\partial  \vec{q}_k} +
    \frac{1}{2} \vec{Q}^-(\vec{q}_k, \vec{q}_{k+1}^n) +
    \frac{1}{2} \vec{Q}^p(\vec{q}_k, \vec{p}_k)
    \right],
  \label{eq:dem:update:q_vec}
\end{equation}
and the stiffness matrix,
\begin{equation}
  \mat{K}(\vec{q}_k, \vec{q}^n_{k+1}, h) =
  -\frac{1}{h}\mat{M}
  -
  \frac{\partial}{\partial \vec{q}_{k+1}} \vec{Q}^-_d(\vec{q}_k, \vec{q}_{k+1}, h).
  \label{eq:dem:update:q_mat}
\end{equation}
The momentum update equation (\ref{eq:update:p}) simplifies to
\begin{equation}
  \vec{p}_{k+1} = D_2L^\alpha_d =
  \frac{1}{h} \mat{M} [\vec{q}_{k+1} - \vec{q}_{k}] - h \alpha \frac{dV}{dq_k}.
  \label{eq:dem:update:p}
\end{equation}
For the first-order integrator ($\alpha = 0$), this can be  interpreted as the product of the
discrete velocity and the mass and is therefore consistent with a discrete time increment.

\section{Numerical Tests}

The integration scheme outlined above is implemented in the Python programming
language \cite{Python, vanRossum1994}.
The complete algorithm is outlined in Algorithm~\ref{alg:vi}.
A Verlet neighbour list~\cite{Verlet1967} efficiently keeps track of potential
contacts and assists in constructing the residual vector and stiffness matrix.
To simplify the implementation and to focus on key features of the algorithm,
the tangential overlap between particles is not
calculated which restricts the following numerical tests to frictionless particles ($\mu=0$).

Walls are implemented using the Hookean contact model
Eqs.~(\ref{eq:dem_contact_normal}-\ref{eq:dem_contact_tan}) by substituting
the position $\vec{r}_j$ with the wall's normal vector and setting
$\vec{v}_j =\vec{0}$.

\begin{algorithm}
  \caption{The first-order ($\alpha=0$) variational integrator algorithm for DEM.}
  \label{alg:vi}
  Read initial state\;
  Assemble $\mat{M}$\;
  $k \leftarrow 0$\;
  \While{$k < K$}{
    Detect Contacts\;

    $\vec{q}^0_{k+1} \leftarrow \vec{q}_k + h\mat{M}^{-1}\vec{p}_{k}$
    \tcc*[r]{The initial guess for the position}
    $n \leftarrow 0$\;
    \While{$E_{n}/E_{n-1} < tol$} {
      $\vec{F} \leftarrow$ Vector term in (\ref{eq:dem:update:q_vec})\;
      $\mat{K} \leftarrow$ Matrix term in (\ref{eq:dem:update:q_mat})\;
      $\Delta \vec{q}^{\Delta n}_{k+1} \leftarrow CG(\mat{K}, -\vec{F})$
      \tcc*[r]{Solve with the conjugate gradient method}
      $\vec{q}^{n+1}_{k+1} \leftarrow \vec{q}^n_{k+1} +  \Delta \vec{q}^{\Delta n}_{k+1}$
      \tcc*[r]{Update the guess for the next iteration}
      $E_{n-1} \leftarrow E_n$\;
      $E_n \leftarrow |\ \Delta \vec{q}^{\Delta n}_{k+1}|$\;
      $n \leftarrow n+1$\;
    }
    $\vec{q}_{k+1} \leftarrow \vec{q}^n_{k+1}$
    \tcc*[r]{Set the coordinates for the next iteration}
    $\vec{p}_{k+1} \leftarrow \frac{1}{h} \mat{M} \left[\vec{q}_{k+1} - \vec{q}_{k}\right]$
    \tcc*[r]{Update the momentum using Eq. (\ref{eq:dem:update:p})}
    $k \leftarrow k+1$\;
  }
\end{algorithm}

\begin{figure}
  \begin{subfigure}[b]{0.33\columnwidth}
    \centering
    \includegraphics[width=0.6\textwidth]{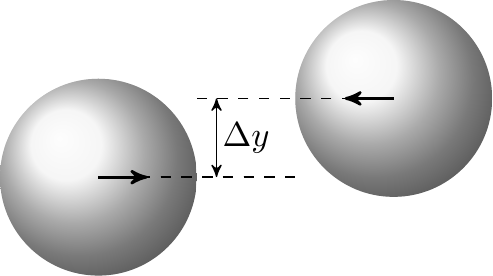}
    \caption{Impact}
    \label{fig:setups:impact}
  \end{subfigure}%
  \begin{subfigure}[b]{0.33\columnwidth}
    \centering
    \includegraphics[width=0.3\textwidth]{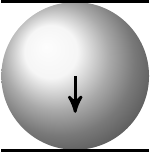}
    \caption{Restitution between walls}
    \label{fig:exp:restitution}
  \end{subfigure}%
  \begin{subfigure}[b]{0.33\columnwidth}
    \centering
    \includegraphics[width=0.6\textwidth]{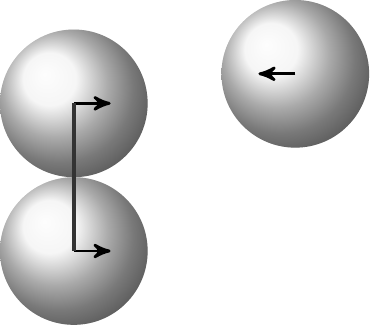}
    \caption{Impact with a bonded pair}
  \end{subfigure}

  \begin{subfigure}[b]{\columnwidth}
    \hfill%
    \includegraphics[width=0.1\textwidth]{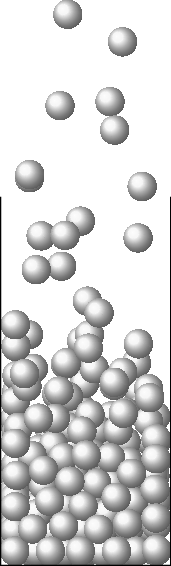}\hspace{1.5cm}%
    \includegraphics[width=0.1\textwidth]{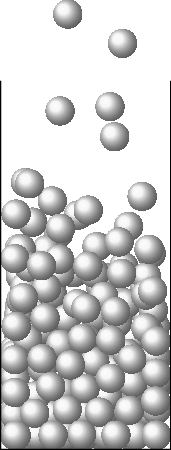}\hspace{1.5cm}%
    \includegraphics[width=0.1\textwidth]{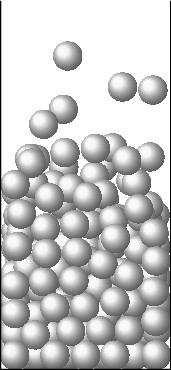}\hfill%
    \hfill
    \caption{Filling a box}
    \label{fig:setups:box}
  \end{subfigure}%
  \caption{
    Particle configurations used for numerical experiments in this section:
    (a) collision between two particles,
    (b) a particle bouncing between walls,
    (c) collision with simplified bonded particles and
    (d) particles filling a box.
    Solid lines between particle centres indicate that a simplified 'bond'
    was present between particles.
  }
  \label{fig:setups}

\end{figure}

A number of numerical experiments are performed to test the integrator and its
implementation.
Figure~\ref{fig:setups} shows the various configurations used in the tests:
a collision between two particles,
a single particle bouncing between two parallel walls,
a collision between a bonded pair and a third particle,
and an ensemble of particles settling in a box under gravity.
Each case is discussed in the following sections.
In each test case all particles have the same diameter and mass and
the same parameters for the Hookean contact model with no friction
($\mu=0$) and
normal and tangential damping is fixed to $\gamma_t = \gamma_n/2$.
A damping parameter is introduced $\gamma = \gamma_n/m_\mathrm{eff}$,
and different values of the parameter $\gamma$ are used to test various aspects of
the integrator.
In all simulations the value of the contact stiffness is fixed relative to other model parameters
such that $kd/mg=\SI{195000}{}$, where $g$ is gravitational acceleration and $m$ and $d$ are
the particle mass and diameter, respectively.

The algorithm outlined above provides an integrator for DEM in a variational
setting. In order to demonstrate that it indeed recovers the same solution as a
conventional DEM simulation,
comparisons are made with the results from the open source software package
LAMMPS~\cite{Plimpton1995}, where possible.
LAMMPS implements a Hookean contact model \cite{Silbert2001c, Brilliantov1996a, Zhang2005}.
The default velocity-Verlet~\cite{Verlet1967} integrator in LAMMPS is used.

\subsection{Two particle impact}

The first test shows the numerical integration of the Hookean contact model over one collision
between two particles.
Simulations for different values of the integration time step $h$, damping
$\gamma$ and offset $\Delta y$ (see Fig~\ref{fig:setups:impact}) are performed.
The two particles have initial positions $\vec{r} \equalhat [d, \pm\Delta y/2, 0]$ and
velocities $\vec{v} \equalhat [\mp v, 0, 0]$. 

An analytical solution is available for the special case when $\Delta y = 0$,
as all contributions from tangential forces remain zero.
The time at which the collision starts can be calculated and is given by $t_A=d/2v$.
After this time, the force between the particles is given by
\begin{equation}
  F = k_n[d - x] +\frac{1}{2} \gamma m v ,
\end{equation}
where $[x, 0, 0]$ and $[-v, 0, 0]$ is the position and velocity, respectively, of the
particle on the right.
This is the same force as a damped simple harmonic oscillator for which the position
and velocity are given by:
\begin{align}
  x(t) &=
  \frac{d}{2} - vt_\gamma \exp\left(-\frac{\gamma\ t}{m}\right)
  \sin\left( \frac{t}{t_\gamma}\right),
  \label{eq:impact:analytic:x}
  \\
  v(t) &= v \exp\left(-\frac{\gamma t }{ m}\right)
  \left[
    \frac{  \gamma^2 }{ m} \sin\left( \frac{t}{t_\gamma} \right)
   -   \cos\left( \frac{t}{t_\gamma} \right)
  \right],
  \label{eq:impact:analytic:v}
\end{align}
where
$t_\gamma = \left[2k/m - [\gamma/m]^2\right]^{-1/2}$.
The duration of the collision can also be calculated by solving for $t$ in $x(t) = d/2$,
which gives $t_C = \pi\sqrt{m/(2k)}$.

\begin{figure}
  \begin{subfigure}[T]{0.33\columnwidth}
    \centering
    \includegraphics{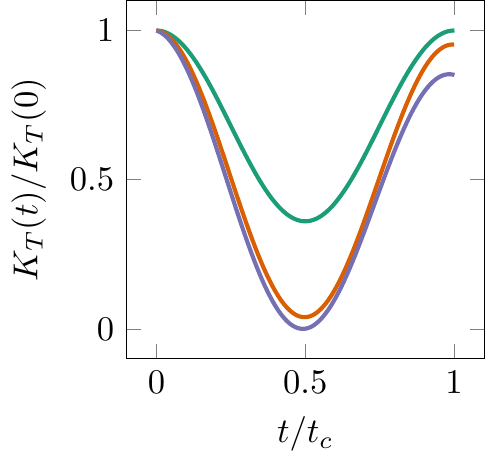}
    \caption{}
    \label{fig:impact:smdt:kl}
  \end{subfigure}
  \begin{subfigure}[T]{0.33\columnwidth}
    \centering
    \includegraphics{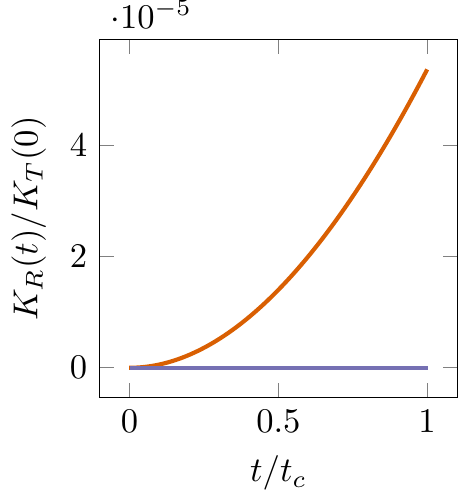}
    \caption{}%
    \label{fig:impact:smdt:kr}
  \end{subfigure}
  \begin{subfigure}[T]{0.33\columnwidth}
    \centering
    \includegraphics{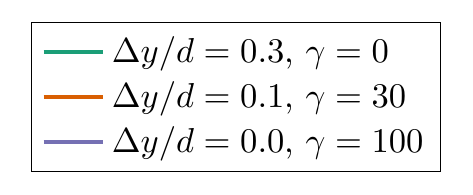}
  \end{subfigure}
  \caption{Two particle impact: Kinetic energy of one particle for small integration time
    steps ($h \approx t_c/160$) and various parameters.
    Results between the proposed integrator and LAMMPS are indistinguishable.
  }
  \label{fig:impact}

\end{figure}

\begin{figure}
  \centering
  \begin{subfigure}[T]{0.33\columnwidth}
    \centering
    \includegraphics{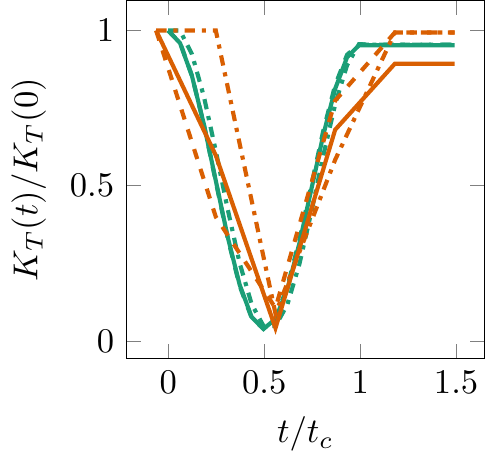}
    \caption{}
    \label{fig:impact:kl}
  \end{subfigure}
  \begin{subfigure}[T]{0.33\columnwidth}
    \centering
    \includegraphics{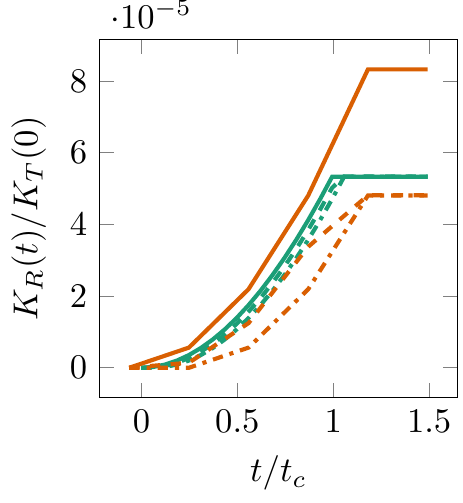}
    \caption{}%
    \label{fig:impact:kr}
  \end{subfigure}
  \begin{subfigure}[T]{0.33\columnwidth}
    \centering
    \includegraphics{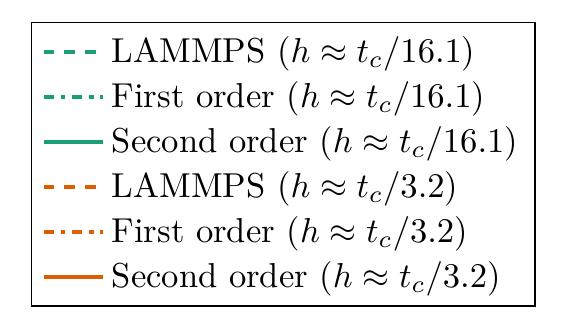}
  \end{subfigure}
  \caption{Two particle impact: Kinetic energy of one particle for large integration time
    steps ($h \approx {t_c/3.2, t_c/16.1}$) and $\Delta y/d = 0.1$ and $\gamma=30$.
    There is agreement between LAMMPS and the second-order integrator for
    $h\approx t_c/16.1$, but the energy dissipation over the course of the collision is not
    correctly calculated for the very large time step $h\approx t_c/3.2$.
  }
  \label{fig:impact:smdt}
\end{figure}

\begin{figure}
  \centering
  \begin{subfigure}[T]{0.33\columnwidth}
    \centering
    \includegraphics{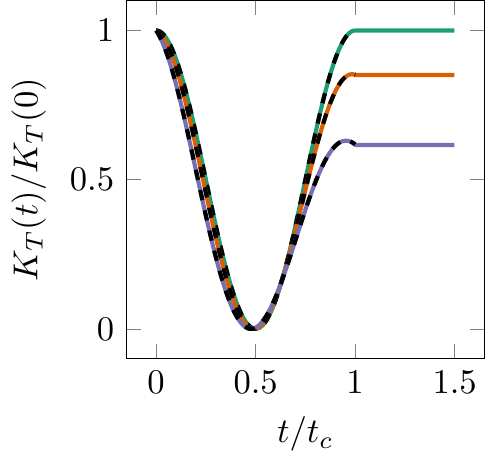}
  \end{subfigure}
  \begin{subfigure}[T]{0.33\columnwidth}
    \centering
    \includegraphics{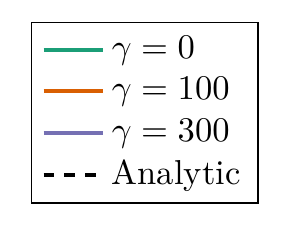}
  \end{subfigure}

  \caption{Two particle impact: Kinetic energy of one particle
    ($h \approx t_c/160$, $\Delta y=0$) compared to the analytic solution in
    equation~(\ref{eq:impact:analytic:v}).
  }
  \label{fig:impact:analytic}
\end{figure}

Figures~\ref{fig:impact} to \ref{fig:impact:analytic} show the translational and rotational
kinetic energy 
of the particle on the right over the course of the collision.
The translational kinetic energy was calculated as $K_T = 1/2\ m[v_x^2+v_y^2]$, where
$\vec{v} \equalhat [v_x, v_y, 0]$ is the particle velocity, and
rotational kinetic energy $K_R = 1/2\ I \omega_z^2$, where $I=2/5\ m [d/2]^2$ is the moment
of inertia of a sphere and $\omega_z$ is the rotational velocity around the $z$ axis.
Figures~\ref{fig:impact:smdt:kl} and \ref{fig:impact:smdt:kr}
shows results of the proposed integrator for a small time step $h\approx t_c/160$.
For this time step size there is excellent agreement with LAMMPS and results between the two
integrators are indistinguishable.
Comparisons between the first and second-order integrators and LAMMPS at larger time steps
($h\approx {t_c/3.2, t_c/16.1}$, $\Delta y/d = 0.1$ and $\gamma=30$)
are made in Figure~\ref{fig:impact:kl} and \ref{fig:impact:kr}.
The second-order integrator compares well with LAMMPS in the case when
$h = t_c/16.1$,
but the amount of energy dissipated is (not surprisingly) incorrectly calculated
for large time steps $h = t_c/3.2$ by all integrators tested.
Finally,
Figure~\ref{fig:impact:analytic} compares the second-order integrator (with $h\approx t_c/160$ and $\Delta y = 0$) to the analytic
solution in equation~(\ref{eq:impact:analytic:v}) and demonstrates near exact agreement.

\subsection{Particle bouncing between walls}

Variational integrators are known to display excellent energy conservation,
despite the fact that energy conservation is not guaranteed \cite{Kane2000}.
To test the energy conservation behaviour of the variational integrator, a simulation
is performed where a particle is placed between two parallel walls set $1.01d$ apart.
The particle's initial velocity is perpendicular to
them (see Figure~\ref{fig:exp:restitution}).
In the undamped case ($\gamma_t = \gamma_n = 0$), the particle will bounce between the
walls without loss of energy and thereby provide a good test for the energy conserving
properties of the integrator.
During the brief periods of no contact the total energy in the system will be the particle's
kinetic energy and during a collision some energy will be converted to potential energy
$V = 1/2\ k_n \delta^2$, where $\delta$ is the overlap between the particle and wall.

\begin{figure}
  \centering
  \begin{subfigure}[T]{0.33\columnwidth}
    \centering
    \includegraphics{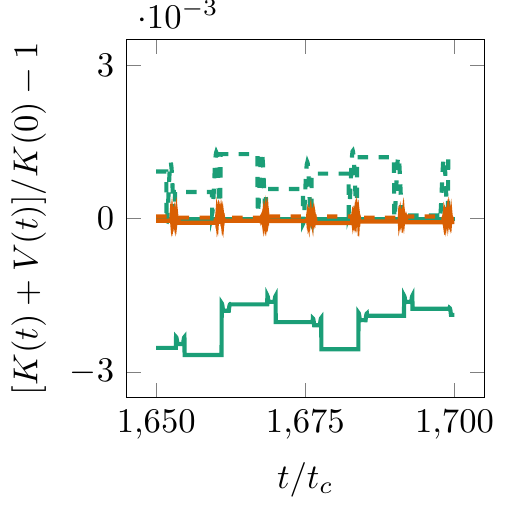}
    \caption{}
    \label{fig:restitution:E}
  \end{subfigure}%
  \begin{subfigure}[T]{0.33\columnwidth}
    \centering
    \includegraphics{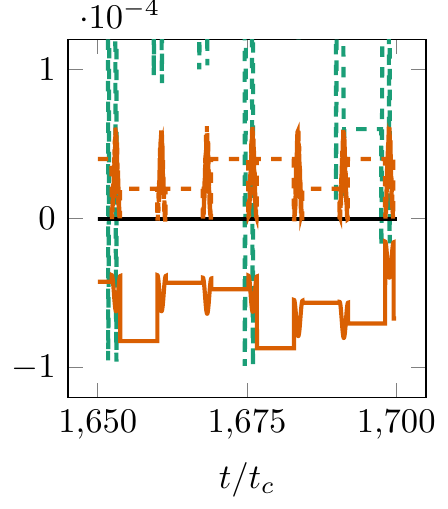}
    \caption{}
    \label{fig:restitution:Ezoom}
  \end{subfigure}%
  \begin{subfigure}[T]{0.33\columnwidth}
    \centering
    \includegraphics{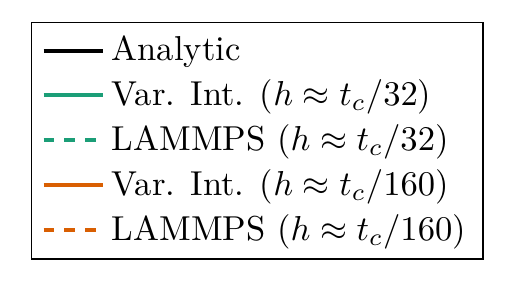}
  \end{subfigure}
  \caption{
    The total energy of a particle bouncing between walls after $250$ collisions.
    A magnified section of (a) is shown in (b).
    Some energy is lost by the variational integrator, but for small time steps
    this is equivalent to the bounded fluctuation produced by the velocity-Verlet
    integrator used in LAMMPS.
  }
  \label{fig:restitution}
\end{figure}

The total energy  of the particle is plotted in Figure~\ref{fig:restitution}.
The total energy is the sum of the kinetic energy and potential energy of the Hookean
contact.
The simulation is carried out over $250$ collisions, but the graph shows the total energy
for the last few collisions.
For $h\approx t_c/32$, the second-order integrator looses a small fraction ($\sim0.1\%$) of energy over the course of the simulation.
The total energy of the LAMMPS simulation is not exactly conserved, but remains bounded.
However, the energy loss of the variational integrator is dependent on the time step
and when the time step is reduced to $h\approx t_c/160$, the energy loss in the
variational integrator is similar to the fluctuations in the energy produced by the
velocity-Verlet integrator used in LAMMPS.

\subsection{Impact with a bonded pair}

A simple bonded particle contact model is implemented by allowing attractive forces between
particles.
This is implemented by creating a 'bond' between particles if in the initial configuration
they are close together ($|\delta_{ij}| < d/100$).
Whenever a bond exists between particles, the potential $V_{ij}=k_n [\delta_{ij}]^2/2$ was
used even when particles were separated.

A simulation is performed of a collision between a pair of bonded particles and a third unbonded particle.
The purpose of this is to test the simplified bonded particle model and
test the integrator with contact models that have different time scales.
Different time scales can be introduced by choosing different spring stiffness constants
for regular Hookean interactions ($k$) and bonded contacts ($k_B$).

\begin{figure}
  \centering
  \begin{subfigure}[T]{0.33\columnwidth}
    \centering
    \includegraphics{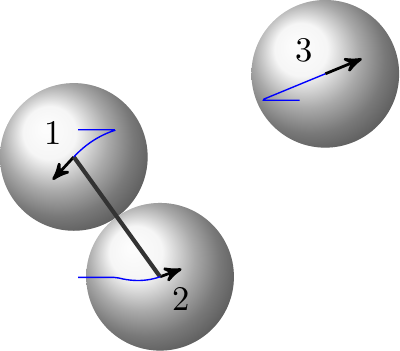}
    \caption{}
    \label{fig:impact:bpm:b}
  \end{subfigure}%
  \begin{subfigure}[T]{0.33\columnwidth}
    \centering
    \includegraphics{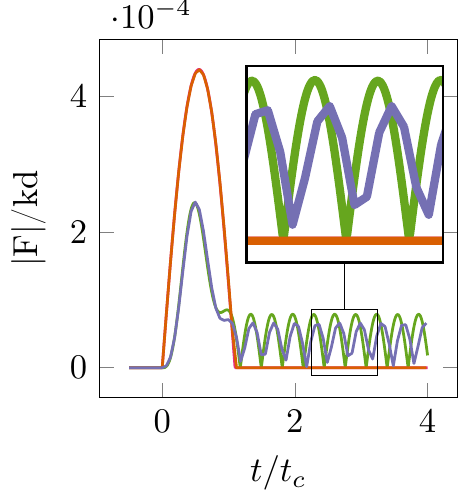}
    \caption{}
    \label{fig:impact:bpm:F}
  \end{subfigure}

  \begin{subfigure}[T]{0.6\columnwidth}
    \centering
    \includegraphics{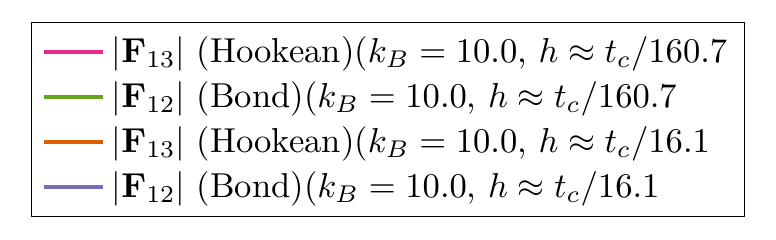}
  \end{subfigure}
  \caption{Impact with a bonded pair: the configuration of particles (a) before and (b)
    after the collision.
    The trajectories of the particles are denoted by  blue lines.
    The inter-particle forces are compared in (c),  particles \SI{1}{} and \SI{2}
    experience a force due to the bond between them and a Hookean contact force
    is present momentarily between \SI{1}{} and \SI{3}{}.}
  \label{fig:impact:bpm}
\end{figure}

The initial setup is similar to the two particle impact simulation, except that one
of the particles is replaced by a bonded pair, see Figure~\ref{fig:impact:bpm}.
The bonded particles are given the same initial velocity.

The bond between the two particles on the left prevents them from separating after
the impact, Figure~\ref{fig:impact:bpm:b} shows the particles and their trajectories after
the collision.
The magnitudes of two inter-particle forces are shown in Figure~\ref{fig:impact:bpm:F}.
The collision between particles $1$ and $2$ produces a peak at the impact.
After the collision, the bond produces an oscillation in the force between particles $1$
and $3$.
The maximum integration time step size is determined by the smallest time scale
($\min\{\sqrt{2k/m}, \sqrt{2k_B/m}\}$).

\subsection{Filling a box}

To investigate a less academic test case, the variational integration scheme is
employed to simulate an ensemble of particles.
A LAMMPS simulation was run to create an initial condition consisting of $N_p=218$ particles
in a $L\times L \times 20L$ (with $L=6d$) box.
A gravitational force was applied in the $[0, 0, -1]$ direction
A snapshot of the simulation captured before all the particles had settled in the bottom of the box was used as the starting configuration of further tests.
The simulation was continued in LAMMPS and the variational integrator until all the
particles settled at the bottom of the box.

When making comparisons between simulations with ensembles of particles, the sensitivity of
these systems to initial conditions and small numerical errors must be kept in mind.
Instead of focusing on individual particle positions and velocities,
macroscopic quantities are compared.
Here, the average kinetic energy per particle and velocity fluctuations
(which is related to the granular temperature) of the ensemble are
shown as the simulation progress.
These quantities are calculated as
\begin{align}
  \bar{K} &= \frac{1}{2N_p}  \sum_{i=1}^{N_P}
  \left\{
    m\left[[v^i_x]^2 + [v^i_y]^2 + [v^i_z]^2\right]
    +
    I\left[\omega^2_x + \omega^2_y + \omega^2_z\right]
  \right\}, \\
  \delta v &=
  \frac{1}{3N_P}
  \sum_{i=1}^{N_P}
  \left[
    [\bar{v}_x - v^i_x]^2 +
    [\bar{v}_y - v^i_y]^2 +
    [\bar{v}_z - v^i_z]^2
    \right],
\end{align}
where bars denote average velocity components: $\bar{v}_x = \sum_i^{N_P} v^i_x/N_P$.

\begin{figure}
  \centering
  \begin{subfigure}[T]{0.33\columnwidth}
    \centering
    \includegraphics{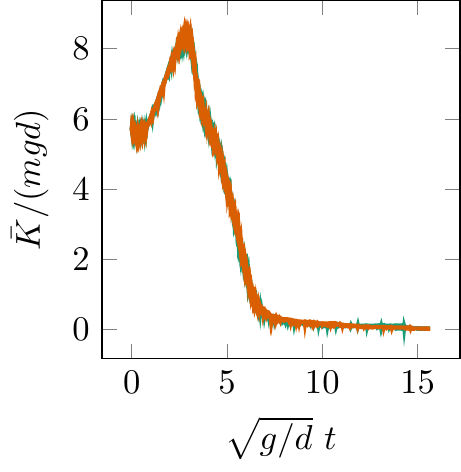}
    \caption{}
    \label{fig:box:K}
  \end{subfigure}
  \begin{subfigure}[T]{0.33\columnwidth}
    \centering
    \includegraphics{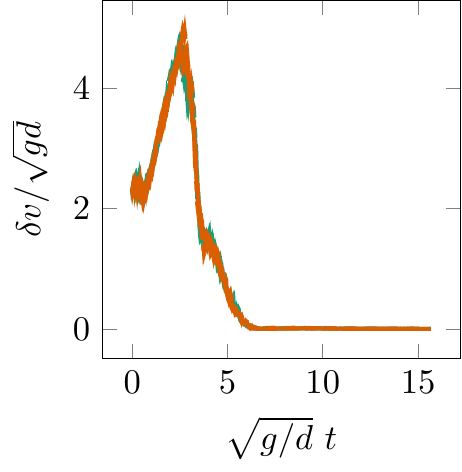}
    \caption{}
    \label{fig:box:T}
  \end{subfigure}%
  \begin{subfigure}[T]{0.33\columnwidth}
    \centering
    \includegraphics{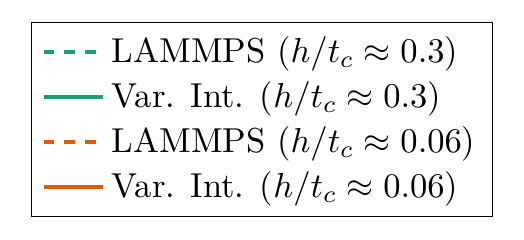}
  \end{subfigure}
  \caption{
    Macroscopic quantities (a) kinetic energy and (b) kinetic stress of $N_p=218$ particles
    settling in a box under gravity.
  }
  \label{fig:box}
\end{figure}

The results are presented in Figure~\ref{fig:box} as a function of the simulation time.
A particle system such as this is known to be sensitive to initial conditions and numerical
errors, so particle trajectories diverge after a few collisions even for the same integration
method with different time steps.
However, the physically meaningful values such as the the coarse-grained statistical quantities
presented show excellent agreement with LAMMPS.

\section{Conclusion}

A variational integrator for DEM has been described and implemented for the Hookean contact model. Our implicit scheme has been compared against the velocity-Verlet method implemented in LAMMPS. Excellent accuracy has been observed at the micro scale (integration over a single collision), macro scale (particles setting in box) as well as good long-term stability (particle bouncing between walls). A simplified bonded particle model has been implemented, thereby demonstrating the method’s versatility and the ability to include other contact models.

Using an implicit numerical method, there is additional computational expense when compared to explicit methods. However, as a variational integrator our approach is attractive since it is a discrete realisation of the Lagrange-d’Alembert principle, an extension of Hamilton’s principle to non-conservative systems, that computes the trajectories of particles by finding the stationary point of the action. Therefore, it represents a dynamical extension of the atomistic
simulations based on the quasi-static energy minimisation principle that inspired the Quasicontinuum (QC) method.
Thus, in a fully realised granular QC method, as motivated in Figure~\ref{fig:outline}, the computational cost of using an implicit integration scheme will be offset by the reduced degrees of freedom of the simulation.
Indeed, in our future work will focus on developing a suitable granular QC method, including appropriate summation rules.

\section*{Acknowledgements}%
This work was supported by the UK Engineering
and Physical Sciences Research Council grant EP/R008531/1 for the
Glasgow Computational Engineering Centre.

\bibliographystyle{model1-num-names}
\bibliography{qc}

\end{document}